**Notes on the Diophantine Equation** $A^4 + aB^4 = C^4 + aD^4$

Paul A. Roediger (proediger@comcast.net)

**Abstract**. A new formulation of the subject equation is presented. Several parametric and semi-parametric solutions are derived. The parametric solution for $a = -1$ was originally presented in 1972, but never published. A computer-generated version was found by Zajta and published in 1983. The pencil-and-paper version is presented for the record, along with many other examples, including eight coincidental triplet solutions.

Let

(1) $$A = p + q,\ C = p - q,\ D = r + s \text{ and } B = r - s.$$

The subject equation becomes

(2) $$a = \frac{pq\left(p^2 + q^2\right)}{rs\left(r^2 + s^2\right)}$$

Let

(3) $$p = ry,\ s = qx \text{ and } r = qt.$$

Equation (2) becomes

(4) $$\frac{ax^3 - y}{y^3 - ax} = t^2$$

Observe that equation (4) may be written as

(5) $$\frac{x - \rho y}{y - a\rho x} = t^2, \text{ where}$$

(6) $$\rho = \frac{xy + 1}{ax^2 + y^2}$$

Solutions of (5) are given by

(7) $$x = \frac{t^2 + \rho}{\omega},\quad y = \frac{a\rho t^2 + 1}{\omega}$$

Substituting (7) into (6) yields:

(8) $$a^2\rho^3 t^4 + \left(3a\rho^2 - 1\right)t^2 + a\rho^3 = \omega^2$$

Solutions of (2) take the form

(9) $$p = t(a\rho t^2 + 1),\ q = \omega,\ r = \omega t \text{ and } s = t^2 + \rho$$

**Case 1**: $a = 1$

With $a = 1$ equation (7) has the obvious solution $\rho = 1$, $\omega = t^2 + 1$, which leads to the trivial solution

$$A = D = (t + 1)\left(t^2 + 1\right), B = C = (t - 1)\left(t^2 + 1\right).$$

Following up this case, take $\rho = 1 + z$, obtaining

$$(t^4+1)z^3 + 3(t^4+t^2+1)z^2 + (t^2+1)^2 z + (t^2+1)^2 = \omega^2 \tag{10}$$

Let $\omega = \frac{3}{2}(t^2+1)z + (t^2+1)$ in equation (10). The constant and linear terms in $z$ cancel, leaving

$(t^4+1)z^3 + 3(t^4+t^2+1)z^2 = \frac{9}{4}(t^2+1)^2 z^2$, which has the solution $z = -\frac{3(t^2-1)^2}{4(t^4+1)}$.

Accordingly, $\rho = 1 + z = \frac{t^4+6t^2+1}{4(t^4+1)}$, $\omega = \frac{(t^2+1)(-t^4+18t^{2-1})}{8(t^4+1)}$ and

$$\begin{aligned} p &= 2t(t^6+10t^4+t^2+4) & q &= (t^2+1)(-t^4+18t^2-1) \\ r &= t(t^2+1)(-t^4+18t^2-1) & s &= 2(4t^6+t^4+10t^2+1) \end{aligned} \tag{11}$$

**Table 1**. Euler's first solution of $A^4+B^4=C^4+D^4$ obtained from equation (11)

| $t$ | $A$ | $B$ | $C$ | $D$ |
|---|---|---|---|---|
| 3 | 158 | -59 | 133 | 134 |
| 2 | 1203 | -76 | 653 | 1176 |
| 5 | 3351 | -2338 | 3494 | 1623 |
| 5/3 | 17332 | 529 | 6673 | 17236 |

**Note**. The free parameter here is $t = \frac{B+D}{A-C}$

Equation (10) is also satisfied by taking $\omega = \frac{3(t^2-1)^2}{8(t^2+1)}z^2 + \frac{3}{2}(t^2+1)z + (t^2+1)$.

In this case, the constant, linear and quadratic terms in terms in $z$ cancel, leaving

$$(t^4+1)z^3 = \frac{9(t^2-1)^4 z^4}{64(t^2+1)^2} + \frac{9(t^2-1)^2 z^3}{8}\text{, so that } z = -\frac{8(t^2+1)^2(-t^4+18t^2-1)}{9(t^2-1)^4}.$$

In this case,

$$\rho = 1 + z = \frac{t^8+92t^6+326t^4+92t^2+1}{9(t^2-1)^4},$$

$$\omega = \frac{(t^2+1)(-t^{12}+214t^6+2481t^8+2804t^6+2481t^4+214t^2-1)}{27(t^2-1)^6}\text{, and}$$

$$\begin{aligned} p &= 3t(t^2-1)^2(t^8+100t^6+190t^4-44t^2+9), \\ q &= (-t^{12}+214t^6+2481t^8+2804t^6+2481t^4+214t^2-1) \\ r &= t(-t^{12}+214t^6+2481t^8+2804t^6+2481t^4+214t^2-1) \\ s &= 3t(t^2-1)^2(9t^8-44t^6+190t^4+100t^2+1) \end{aligned} \tag{12}$$

**Table 2**. Euler's second solution of $A^4+B^4=C^4+D^4$ obtained from equation (12)

| $t$ | $A$ | $B$ | $C$ | $D$ |
|---|---|---|---|---|
| 3 | 10381 | 10203 | 2903 | 12231 |
| 2 | 1584749 | 2061283 | -555617 | 2219449 |
| 5 | 2533177 | 1123601 | 2367869 | 1834883 |

**Note**. The free parameter again is $t=\dfrac{B+D}{A-C}$

Lander (1968) provided methods to derive new parametric solutions of $A^4+B^4=C^4+D^4$ when one starts out with a known solution. For example, starting with (11), Euler's first solution, he derived Euler's second solution (equation (12)), and another solution of degree nineteen. Brudno (1969) obtained two more solutions of degree nineteen and thirty-one.

**Case 2**: $a=-1$

With $a=-1$, solving equation (8) for $t^2$ eventually leads to other parametric solutions, as follows:

(13) $$t^2=\frac{3\rho^2+1+\Delta}{2\rho^3}\text{, where}$$

(14) $$\Delta^2=\left(\rho^2+1\right)^2\left(4\rho^2+1\right)+4\rho^3\omega^2$$

With $\omega=\dfrac{\rho^2+1}{\rho}z$, equation (14) becomes

(15) $$\frac{\Delta^2}{\left(\rho^2+1\right)^2}=4\rho^2+4\rho z^2+1$$

Since the right side of equation (15) is a perfect square when $z=1$, let $z=1+k$.

A rational $k$ is obtained by setting $\dfrac{\Delta}{\left(\rho^2+1\right)}=2nk-\left(2\rho+1\right)$, namely $k=\dfrac{\left(2\rho+1\right)n+2}{\rho n^2-1}$.

Accordingly, $$z=\frac{\rho n^2+\left(2\rho+1\right)n+1}{n^2\rho-1},$$

$$\omega=\frac{(\rho^2+1)(\rho n^2+\left(2\rho+1\right)n+1)}{p(n^2\rho-1)},$$

$$\Delta=\frac{\left(\rho^2+1\right)\left(\left(2\rho+1\right)\rho n^2+4\rho n+\left(2\rho+1\right)\right)}{\rho n^2-1}\text{, and}$$

(16) $$t^2=\frac{n^2\rho^3+\left(2n^2+2n+1\right)\rho^2+\left(n^2-1\right)\rho+\left(n+1\right)^2}{\rho^2\left(n^2\rho-1\right)}$$

Let $t=\dfrac{\rho+v}{\rho}$. The $\rho^3$ term in equation (16) drops out, as does the $\rho^2$ term if $v=\dfrac{n^2+n+1}{n^2}$.

The remaining linear equation has solution

$$\rho = \frac{v^2 + (n+1)^2}{n^2 v^2 - 2v - (n^2 - 1)} = \frac{n^6 + 2n^5 + 2n^4 + 2n^3 + 3n^2 + 2n + 1}{n^2 (2n^3 + 2n^2 - 1)}$$ , so that

$$t = \frac{n(n^5 + 4n^4 + 6n^3 + 6n^2 + 4n + 1)}{n^6 + 2n^5 + 2n^4 + 2n^3 + 3n^2 + 2n + 1}$$ , and

$$\omega = \frac{(n^5 + 2n^4 + 2n^3 + 2n^2 + n + 1)(n^8 + 4n^7 + 12n^6 + 20n^5 + 21n^4 + 16n^3 + 10n^2 + 4n + 1)}{n^3 (2n^3 + 2n^2 - 1)(n^6 + 2n^5 + 2n^4 + 2n^3 + 3n^2 + 2n + 1)}$$ . and

(17)

$$\begin{aligned}
p &= -n^4 (n+1)(n^2 + 2n + 2)(n^4 + 3n^3 + 3n^2 + 3n + 1) \\
q &= (n^2 + n + 1)(n^3 + n^2 + 1)(n^6 + 2n^5 + 2n^4 + 2n^3 + 3n^2 + 2n + 1) \\
r &= n(n+1)(n^2 + n + 1)(n^3 + n^2 + 1)(n^4 + 3n^3 + 3n^2 + 3n + 1) \\
s &= n (n^{10} + 4n^9 + 8n^8 + 10n^7 + 7n^6 + 2n^5 + n^4 + 2n^3 + 3n^2 + 2n + 1)
\end{aligned}$$

Zajta (1983) surveys known transformation methods leading to solutions of $A^4 + B^4 = C^4 + D^4$ . There several new methods are documented, including so-called Simple-Dual (SD) and Composite-Dual (CD) parametric transformations. By applying an *unspecified* CD transformation to a variant of equation (11), the simplest of known parametric solutions, he presents $P_1(u,v)$ , which is a version of equation (17) with $n = -\dfrac{u+v}{v}$ .

**Table 3**. Solutions of $A^4 - B^4 = C^4 - D^4$ obtained from equation (17)

| $n$ | $A$ | $B$ | $C$ | $D$ |
|---:|---:|---:|---:|---:|
| 1 | 7 | 157 | -227 | 239 |
| -2 | -257 | 292 | 193 | -256 |
| -1/2 | 502 | 298 | -497 | -271 |
| -3/2 | -6842 | 9018 | -4903 | -8409 |
| 1/2 | 6742 | 5098 | -9043 | 8531 |
| 2 | -10757 | 18292 | -45883 | 46136 |
| -3 | -28997 | 33237 | 59777 | -60369 |
| -1/3 | 89841 | 27879 | -90829 | -43307 |

**Note**. The free parameter is $n = \dfrac{(y-x)(y+1)}{x-(y^2+y+1)}$ , where $x = \dfrac{D-B}{A-C}$ and $y = \dfrac{A+C}{D+B}$

Following the method of Lander (1968), two more solutions of the $a = 1$ case may be derived from equation (17). They are of degree thirteen and fifteen, respectively as follows:

(18)

$$\begin{aligned}
p &= (n^3 + n^2 + 1)(n^2 + n + 1)(n^8 + 4n^7 + 9n^6 + 14n^5 + 14n^4 + 10n^3 + 6n^2 + 2n + 1) \\
q &= n^4 (n^2 + 2n + 2)(n^3 + n - 1)(n^4 + 2n^3 + 2n^2 + n + 1) \\
r &= n (n^{12} + 6n^{11} + 19n^{10} + 40n^9 + 64n^8 + 80n^7 + 82n^6 + 68n^5 + 46n^4 + 26n^3 + 12n^2 + 4n + 1) \\
s &= n(n^3 + n^2 + 1)(n^3 + n - 1)(n^2 + n + 1)(n^4 + 2n^3 + 2n^2 + n + 1)
\end{aligned}$$

$$p=-(n+1)\begin{pmatrix} n^{14}+8n^{13}+32n^{12}+90n^{11}+195n^{10}+320n^{9}+391n^{8}+ \\ 358n^{7}+254n^{6}+146n^{5}+71n^{4}+30n^{3}+12n^{2}+4n+1 \end{pmatrix}$$

(19)
$$q=n(n+1)(n^4+3n^3+3n^2+3n+1)(n^4+2n^3-n-1)(n^5+5n^4+8n^3+5n^2+n+1)$$

$$r=n(n+1)^4(n^4+3n^3+3n^2+3n+1)(n^6+4n^5+9n^4+6n^3+3n^2+2n+1)$$

$$s=(n^4+2n^3-n-1)(n^5+5n^4+8n^3+5n^2+n+1)(n^6+2n^5+2n^4+2n^3+3n^2+2n+1)$$

**Table 4**. More solutions of $A^4+B^4=C^4+D^4$

| Equation | $n$ | $A$ | $B$ | $C$ | $D$ |
|---|---|---|---|---|---|
| (18) | 1 | 292 | 193 | 257 | 256 |
| | -2 | -2797 | 248 | 2131 | -2524 |
| | -1/2 | 2345 | -2986 | 3190 | 1577 |
| | 1/2 | 60763 | 38078 | 62206 | 29531 |
| (19) | -2 | -239 | 7 | -227 | 157 |
| | 1 | 4288 | 4303 | 3364 | 4849 |
| | -1/2 | 2707 | 6730 | 3070 | -6701 |
| | -3/2 | -73703 | 154522 | -151394 | -92839 |

**Note**. $P_2(u,v)$ of Zajta (1983, pp. 651), becomes equation (18) if one sets $n=v/u$

**Case 3**: Parametric solutions of equation (8) with $\rho=1$

Parametric solutions of equation (2) abound in the literature. The most remarkable ones are parametric solutions for a constant value of $a$, where $p$, $q$, $r$ and $s$ depend on a free parameter, say $u$. Many easier to find semi parametric solutions where $a$ is a function of fractional values of $u$, called here $a_f(u)$ or simply $a_f$. An example of this type was noted by Hayashi (1911):

(20)
$$p=u(u^2-3),\ q=2(u^2-1),\ r=u(u^2-1),\ s=2,\ a_f=u^2-3.$$

For $u=\dfrac{7}{4}$ one calculates $a_f=\dfrac{1}{16}$, which leads to an $a=1$ solution discovered much earlier by Euler, namely

$$542^4+103^4=514^4+359^4$$

Many more examples may be found by setting $\rho=1$ and $\omega=at^2-\alpha$ in equation (8), yielding the following version of equation (9):

(21)
$$p=t(at^2+1),\ q=at^2-\alpha,\ r=t(at^2-\alpha),\ s=t^2+1 \text{ and } a_f=\frac{\alpha^2+t^2}{(2\alpha+3)t^2+1}.$$

The $\rho=1$ restriction is not as limiting as first thought: examination reveals that $(ac^{-4},\rho c^2,tc,\omega c)$ satisfies (8) whenever $(a,\rho,t,\omega)$ does. Thus, setting $\rho=1$ is really equivalent to restricting $\rho$ to be a perfect square.

Judicious play with $\alpha$ and $t$ (mostly as functions of $u$) leads to many other examples. Fractional values of $u$ which provide either a constant value for $a_f$ or values of $a_f$ capable of having large quartic factors (leading to relatively small values of $a$) are listed in Table 5. The corresponding solutions are then elaborated upon in Tables 6 through Table 9.

**Table 5**. For $\rho = 1$, some additional $a_f$'s for various choices of $\alpha$ and $t$

| $i$ | $\alpha$ | $t$ | $a_f$ | $i$ | $\alpha$ | $t$ | $a_f$ |
|---|---|---|---|---|---|---|---|
| 1 | $\frac{1}{2}$ | $u$ | $\frac{1}{4}$ | 7 | $\frac{-2}{u^2}$ | $\frac{1}{u}$ | $\frac{1}{u^2-1}$ |
| 2 | $\frac{3u^2+4}{u^2}$ | $u$ | $\frac{u^2+4}{9u^4}$ | 8 | $\frac{u^4+2u^2+2}{2\left(1-u^2\right)}$ | $\frac{3u^2+2}{2u\left(u^2-1\right)}$ | $\frac{1-u^4}{5}$ |
| 3 | $\frac{1}{u^2}$ | $\frac{1}{u}$ | $\frac{1}{u^2+2}$ | 9 | $\frac{9u^2+16}{1-6u^2}$ | $\frac{4u^2+1}{8\left(2-u^2\right)}$ | $\frac{u^2-7}{u^2-3}$ |
| 4 | $-\left(3u^2+4\right)$ | $u$ | $\frac{9u^2+16}{1-6u^2}$ | 10 | $\frac{3}{2}$ | $u$ | $u^2+\frac{9}{4}$ |
| 5 | $\frac{3u^2+4}{u^2\left(u^2+2\right)}$ | $\frac{u}{u^2+2}$ | $\frac{u^2+2}{u^4}$ | 11 | $\frac{2\left(4u^2+1\right)+9u\sqrt{2}}{2\left(2-u^2\right)}$ | $\sqrt{2}$ | $\frac{4u^2+1}{2\left(2-u^2\right)}$ |
| 6 | $\frac{1-4u^2}{4}$ | $u$ | $\frac{4u^2+1}{8\left(2-u^2\right)}$ | | | | |

**Table 6**. Sample $p, q, r, s$ and $a$ from Table 5 (numerical examples provided in Tables 8 and 9)

| $i$ | $p$ | $q$ | $r$ | $s$ | $a_f$ |
|---|---|---|---|---|---|
| 1 | $u\left(u^2+4\right)$ | $u^2-2$ | $u\left(u^2-2\right)$ | $4(u^2+1)$ | $\frac{1}{4}$ |
| 2 | $u\left(u^2+16\right)$ | $u^2-20$ | $u\left(u^2-20\right)$ | $9u^2$ | $(u^2+4)^2/(9u^4)$ |
| 3 | $u^2+1$ | $u$ | $1$ | $u\left(u^2+2\right)$ | $(u^2+2)^{-1}$ |
| 4 | $u\left(9u^2+1\right)$ | $9u^2-4$ | $u\left(9u^2-4\right)$ | $1-6u^2$ | $\frac{9u^2+16}{1-6u^2}$ |
| 5 | $u\left(u^2+1\right)$ | $3\left(u^2+2\right)$ | $3u$ | $u^2\left(u^2+4\right)$ | $(u^2+2)u^{-4}$ |
| 6 | $u\begin{pmatrix}4u^4- \\ 7u^2+16\end{pmatrix}$ | $\begin{pmatrix}4u^4- \\ 19u^2+4\end{pmatrix}$ | $u\begin{pmatrix}4u^4- \\ 19u^2+4\end{pmatrix}$ | $8\left(u^2+1\right)\left(2-u^2\right)$ | $\frac{4u^2+1}{8(2-u^2)}$ |
| 7 | $u^4-u^2+1$ | $u\left(2u^2-1\right)$ | $2u^2-1$ | $u\left(u^4-1\right)$ | $(u^2-1)^{-1}$ |
| 8 | $\left(3u^2+2\right)\left(9u^2+1\right)$ | $2u\left(u^2-1\right)^2$ | $\left(3u^2+2\right)\left(1-u^2\right)$ | $10u\left(4u^2+1\right)$ | $(1-u^4)/5$ |
| 9 | $(3u^2-5)(u^6-u^4-9u^2+25)$ | $u(u^2-7)(u^6-3u^4+3u^2-25)$ | $(3u^2-5)(u^6-3u^4+3u^2-25)$ | $u(u^2-3)(u^2+1)\left(u^4-6u^2+25\right)$ | $\dfrac{u^2-7}{u^2-3}$ |
| 10 | $u\left(4u^4+9u^2+4\right)$ | $4u^4+9u^2+6$ | $u\left(4u^4+9u^2+6\right)$ | $4\left(u^2+1\right)$ | $u^2+\dfrac{9}{4}$ |
| 11 | $2\left(u^2+1\right)\sqrt{2}$ | $3u\sqrt{2}$ | $6u$ | $2\left(2-u^2\right)$ | $\dfrac{4u^2+1}{2(2-u^2)}$ |
| 12 | $2\left(u^2+1\right)$ | $3u$ | $6u$ | $2\left(2-u^2\right)$ | $\dfrac{4u^2+1}{8(2-u^2)}$ |
| 13 | $(u^2-1)(u^4+14u^2+1)$ | $u(u^4-22u^2+1)$ | $(u^2+1)(u^4-22u^2+1)$ | $9u(u^4-1)$ | $\dfrac{1}{9}$ |

**Note**. Line 12 is a rational version of line 11 of Table 5 and line 13 is a special case of line 2 of Table 5

**Table 7**. Numerical examples obtained from Table 6, line 1, satisfying $A^4+aB^4=C^4+aD^4=n\le 10^{12}$

| $a$ | $A$ | $B$ | $C$ | $D$ | $u$ | $a_f$ | $n$ |
|---|---|---|---|---|---|---|---|
| 4 | 7 | 6 | -9 | 4 | 1/1 | 1/4 | 7585 |
| | 33 | 31 | -47 | 3 | 1/2 | 1/4 | 4880005 |
| | 61 | 32 | 19 | -46 | 2/3 | 1/4 | 18040145 |
| | 103 | 88 | -137 | -14 | 1/3 | 1/4 | 352429025 |
| | 219 | 122 | -11 | -168 | 2/5 | 1/4 | 3186391345 |
| | 241 | 189 | -303 | -59 | 1/4 | 1/4 | 8477361925 |
| | 317 | 171 | 147 | -239 | 4/5 | 1/4 | 13518183445 |
| | 331 | 311 | -469 | 127 | 3/4 | 1/4 | 49423420085 |
| | 471 | 346 | -569 | -144 | 1/5 | 1/4 | 106541111905 |
| | 529 | 324 | -129 | -418 | 2/7 | 1/4 | 122390827585 |
| | 711 | 373 | 137 | -537 | 4/7 | 1/4 | 332978996005 |
| | 557 | 532 | -803 | 122 | 3/5 | 1/4 | 416664780305 |
| | 817 | 571 | -959 | -281 | 1/6 | 1/4 | 870752499845 |
| | 913 | 502 | 479 | -688 | 6/7 | 1/4 | 948861341825 |

**Note**. Solutions with $v=i/j$ and $v=j/2i$ for $i<j$ are the same

**Table 8**. Numerical examples obtained from Table 6, line 2, satisfying $A^4+aB^4=C^4+aD^4=n\le 10^{12}$

| $a$ | $A$ | $B$ | $C$ | $D$ | $u$ | $a_f$ | $n$ |
|---|---|---|---|---|---|---|---|
| 9 | -85 | 361 | -625 | 77 | 2/1 | 1/9 | 152904267994 |
| | 250 | 503 | -830 | 329 | 3/1 | 1/9 | 580028236729 |
| | -1105 | 1435 | -2509 | -233 | 3/2 | 1/9 | 39654570456250 |
| | 1513 | 1519 | -2159 | 1367 | 4/1 | 1/9 | 53155527769450 |
| | 7514 | 5761 | -6526 | 5951 | 5/1 | 1/9 | 13101413507132185 |
| | 17723 | 12001 | -10249 | 13213 | 6/1 | 1/9 | 285347996556595450 |

**Note**. Solutions with $u=i/j$ and $u=j/i$ are the same

**Table 9**. Other numerical examples obtained from Table 6 satisfying $A^4 + aB^4 = C^4 + aD^4 = n$

| $a$ | $A$ | $B$ | $C$ | $D$ | $i$ | $u$ | $a_f$ | $n$ |
|---|---|---|---|---|---|---|---|---|
| 1 | 558 | -503 | 222 | 631 | 3 | 7/4 | 16/81 | 160961094577 |
| | 631 | -222 | 503 | 558 | 8 | 1/3 | 16/81 | 160961094577 |
| | 1381 | 878 | 997 | -1342 | 8 | 3 | -16 | 4231525221377 |
| | 2854 | -1797 | -1034 | 2949 | 5 | 7/4 | 1296/2401 | 76773963505537 |
| | 10943964 | -1733885 | -5558948 | 10758915 | 10 | 7/16 | 625/256 | 14353974451590106715110524241 |
| 2 | 248 | -223 | 44 | 257 | 7 | 3 | 1/8 | 6255715457 |
| | -16727 | 36384 | 41513 | 23532 | 7 | 7/9 | -81/32 | 3583152677936711713 |
| 3 | -2 | 3 | 4 | 1 | 3 | 1 | 1/3 | 259 |
| | 11 | -2 | -7 | 8 | 5 | 1 | 3 | 14689 |
| | 7 | 8 | -11 | -2 | 9 | 1 | 3 | 14689 |
| | -23 | 27 | 37 | -1 | 7 | 2 | 1/3 | 1874164 |
| | 93 | -134 | 63 | 136 | 3 | 5 | 1/27 | 397223137 |
| | 277 | 149 | 241 | -191 | 8 | 2 | -3 | 7365992644 |
| | 241 | 191 | 277 | -149 | 9 | 2 | -3 | 7365992644 |
| | 304 | -127 | 268 | 193 | 8 | 1/2 | 3/16 | 9321150979 |
| | -426 | 211 | 444 | 49 | 5 | 5 | 27/625 | 38879896899 |
| | 86 | 997 | 1256 | 631 | 9 | 3 | 1/3 | 2964216377059 |
| | -16703 | 6064 | 16897 | 3348 | 7 | 7 | 1/48 | 81892097850563329 |
| 5 | 3 | 0 | 1 | 2 | 4 | 1 | -5 | 81 |
| | 22 | 17 | 4 | 19 | 11 | 3/2 | -5 | 651861 |
| | 49 | -137 | 197 | -85 | 6 | 3/2 | -5 | 1767141606 |
| | -58879 | 15860 | 59201 | 10064 | 7 | 9 | 1/80 | 12334623789892762881 |
| | -64151 | 34620 | -51031 | 43152 | 7 | 1/9 | -81/80 | 24118655609148141601 |

**Note**. Solutions with $v = i/j$ and $v = j/2i$ for $i < j$ are the same

**Case 5**: The semi-parametric solutions found in Table 7 for $i = 11$ and $i = 12$ may be further generalized. In equation (21), substitute $\alpha = x + y\sqrt{k}$ and $t = \sqrt{k}$ into $a_f = \dfrac{\alpha^2 + t^2}{(2\alpha + 3)t^2 + 1}$ to get $a_f = \dfrac{x^2 + (y^2+1)k + 2xy\sqrt{k}}{(2x+3)k + 1 + 2yk\sqrt{k}}$. Therefore, $a_f = \dfrac{x}{k}$ when $\dfrac{x^2 + (y^2+1)k}{(2x+3)k+1} = \dfrac{2xy\sqrt{k}}{2yk\sqrt{k}}$, i.e., when

$$(ky)^2 = kx^2 + (3k+1)x - k^2 \tag{22}$$

With $k = \dfrac{u}{(u-1)^2}$ the right-hand side of equation (22) has rational factors

$$(uy)^2 = \left((u-1)x + u^2\right)\left(u(u-1)x - 1\right) \tag{23}$$

which admit, by setting $uy = \left(u(u-1)x - 1\right)v$, the following parameterized quantities:

$$\begin{aligned} x &= \frac{u^2+v^2}{(u-1)(uv^2-1)} & \qquad y &= \frac{\left(u(u-1)x-1\right)v}{u} = \frac{v(u^3+1)}{u(uv^2-1)} \\ \alpha &= \frac{u(u^2+v^2) + v(u^3+1)\sqrt{u}}{u(u-1)(uv^2-1)} & \qquad a_f &= \frac{(u-1)(u^2+v^2)}{u(uv^2-1)} \end{aligned} \tag{24}$$

From equation (9), and after some simplification, the following rational solution of (2) is provided:

(25)
$$
\begin{aligned}
&p = u(u^2 - u + 1)(v^2 + 1) \qquad && q = (u-1)(u^3+1)v \\
&r = (u^3+1)v && s = (u^2 - u + 1)(uv^2 - 1) \\
&a_f = \frac{u(u-1)(u^2+v^2)}{(uv^2-1)}
\end{aligned}
$$

Natural number solutions for $a < 10$ are tabulated below.

**Table 10**. Additional numerical solutions of $A^4 + aB^4 = C^4 + aD^4$ (for small $a$ )

| $i$ | $a$ | $A$ | $B$ | $C$ | $D$ | $u$ | $v$ | $a_f$ | $n$ |
|---|---|---|---|---|---|---|---|---|---|
| 1 | 1 | 239 | -7 | -157 | 227 | 10/1 | 5/4 | 625 | 3262811042 |
| 2 | | -193 | -292 | -256 | 257 | -5/4 | 20/7 | -625/256 | 8657437697 |
| 3 | 2 | 61 | 116 | 139 | -34 | 4/1 | 4/5 | 128 | 375973713 |
| 4 | | 257 | 121 | 263 | 43 | 8/5 | 4/5 | 128 | 4791188163 |
| 5 | | -44 | 257 | 248 | 223 | 9/8 | 3/1 | 81/512 | 8728688898 |
| 6 | | -16727 | -36384 | -41513 | 23532 | -49/32 | 7/9 | -194481/32768 | 3583152677936711713 |
| 7 | 3 | 4 | 1 | 2 | 3 | 3/2 | 1/2 | -3 | 259 |
| 8 | | -11 | -2 | 7 | -8 | -2/1 | 2/1 | -16/3 | 14689 |
| 9 | | 23 | 27 | 37 | -1 | 10/3 | 5/12 | -625/3 | 1874164 |
| 10 | | 37 | 1 | 23 | 27 | 4/3 | 2/1 | 16/27 | 1874164 |
| 11 | | -277 | -149 | 241 | 191 | -4/3 | 2/9 | -16/3 | 7365992644 |
| 12 | | 268 | 193 | 304 | 127 | 16/13 | 11/13 | -16/3 | 9321150979 |
| 13 | | 178 | -399 | -526 | 113 | -13/3 | 13/11 | -28561/432 | 77038751059 |
| 14 | | 16897 | -3348 | 16703 | 6064 | 49/48 | 7/1 | 2401/110592 | 81892097850563329 |
| 15 | 5 | 22 | 17 | 4 | 19 | 2/1 | 2/3 | -80 | 651861 |
| 16 | | -71 | 101 | 147 | 63 | 12/7 | 14/3 | 2000/2401 | 545713686 |
| 17 | | 51031 | -43152 | -64151 | 34620 | -80/1 | 9/1 | -6480 | 24118655609148141601 |

**Note**. Equation (24) has the symmetry for $a$ : $u \to u^{-1}$, $v \to v^{-1}$. From equation (25),
$\dfrac{s(u,v)}{r(u,v)} = sign(uv^2 - 1)\dfrac{r(-u,v)}{s(-u,v)}$. Thus, the arguments may be confined to $|u| > 1$, $v \ge 1$ and $u \ne v^{-2}$.

**Case 6**: Coincidental triplet solutions of $A^4 + aB^4 = C^4 + aD^4 = E^4 + aF^4$

Take

(26)
$M = j^2 + ak^2$, $I = j^4 - 6aj^2k^2 + a^2k^4$ and $J = 4jk(j^2 - ak^2)$, so that $M^4 = I^2 + aJ^2$ and
$M^4\left(A^4 + aB^4\right) = E^4 + aF^4$, where $E^2 = A^2I + aB^2J$ and $F^2 = B^2I - A^2J$

A search over known numerical examples of $A^4 + aB^4 = C^4 + aD^4$ reveals eight non-trivial solutions. The details are tabulated below:

**Table 11**. Eight solutions satisfying $M^4(A^4+aB^4)=M^4(C^4+aD^4)=E^4+aF^4$

| Equation (25) | | | | | | | | Equation (26) | | | | |
|---|---|---|---|---|---|---|---|---|---|---|---|---|
| $i$ | $u$ | $v$ | $a$ | $A$ | $B$ | $C$ | $D$ | $j$ | $k$ | $M$ | $E$ | $F$ |
| 1 | -2/1 | 2/1 | 3 | 7 | 8 | 11 | 2 | 4 | 1 | 19 | 197 | 108 |
| 2 | 3/2 | 1/2 | 3 | 2 | 3 | 4 | 1 | 8 | 1 | 67 | 254 | 137 |
| 3 | 12/5 | 8/15 | -4 | 47 | 31 | 33 | 3 | -11 | 4 | 57 | 16407 | 11601 |
| 4 | -2/1 | 1/1 | -10 | 7 | 4 | 1 | 2 | 8 | 1 | 54 | 122 | 112 |
| 5 | 7/6 | 3/2 | 35 | 8 | 1 | 6 | 3 | -20 | 1 | 435 | 1980 | 1395 |
| 6 | -15/2 | 5/12 | -40 | 94 | 37 | 42 | 3 | 35 | 2 | 10065 | 51390 | 16515 |
| 7 | 7/4 | 1/14 | -65 | 361 | 89 | 686 | 238 | -3 | 1 | 56 | 38416 | 13328 |
| 8 | 4/3 | 5/4 | -111 | 79 | 24 | 44 | 11 | 12 | 1 | 33 | 4719 | 1452 |

**Note**: Table represents results of a search for solutions in natural numbers $a$, $|a|<126$

Interestingly, all coincident solutions found started with an $A$, $B$ and $a$ generated by the Case 5 methodology.

Two links (Wroblewski, dates unknown) are provided as references. There, another very interesting method is described for finding "magical" solutions for the $a=-1$ case. This approach appears to be unrelated to the one described here.

**Appendix**. The smallest solutions of $n = A^4 + aB^4 = C^4 + aD^4$ for $a = 2, 3, 4, 5, 9$

**Table A1**. $a = 2$ (11 solutions, $n \le 10^{13}$)

| $i$ | $A$ | $B$ | $C$ | $D$ | $n$ |
|---|---|---|---|---|---|
| 1 | 139 | 34 | 61 | 116 | 375973713 |
| 2 | 141 | 165 | 37 | 175 | 1877655411 |
| 3 | 263 | 43 | 257 | 121 | 4791188163 |
| 4 | 248 | 223 | 44 | 257 | 8728688898 |
| 5 | 309 | 24 | 301 | 146 | 9117284913 |
| 6 | 577 | 662 | 489 | 684 | 494957326113 |
| 7 | 969 | 167 | 849 | 653 | 883203352163 |
| 8 | 984 | 189 | 604 | 797 | 940071661218 |
| 9 | 999 | 312 | 103 | 844 | 1014957704673 |
| 10 | 1414 | 499 | 354 | 1197 | 4121587360818 |
| 11 | 1751 | 236 | 641 | 1466 | 9406566220833 |
|  | 5028 | 1587 | 1376 | 4243 | 651804419270178 |
|  | 10036 | 1777 | 9072 | 6417 | 10164721972794498 |
|  | 41513 | 23532 | 16727 | 36384 | 3583152677936711713 |

**Table A2**. $a = 3$ (23 solutions, $n \le 10^{12}$)

| $i$ | $A$ | $B$ | $C$ | $D$ | $n$ |
|---|---|---|---|---|---|
| 1 | 4 | 1 | 2 | 3 | 259 |
| 2 | 11 | 2 | 7 | 8 | 14689 |
| 3 | 37 | 1 | 23 | 7 | 1874164 |
| 4 | 93 | 134 | 63 | 136 | 1042059009 |
| 5 | 173 | 97 | 163 | 111 | 1161332884 |
| 6 | 209 | 38 | 197 | 108 | 1914285169 |
| 7 | 197 | 108 | 133 | 152 | 1914285169 |
| 8 | 233 | 131 | 41 | 189 | 3830795284 |
| 9 | 263 | 27 | 257 | 109 | 4785944884 |
| 10 | 268 | 67 | 254 | 137 | 5219140339 |
| 11 | 254 | 137 | 134 | 201 | 5219140339 |
| 12 | 277 | 149 | 241 | 191 | 7365992644 |
| 13 | 304 | 127 | 268 | 193 | 9321150979 |
| 14 | 307 | 233 | 89 | 277 | 17724760564 |
| 15 | 356 | 163 | 38 | 279 | 18179748979 |
| 16 | 383 | 110 | 331 | 240 | 21956892721 |
| 17 | 444 | 49 | 426 | 211 | 38879896899 |
| 18 | 526 | 113 | 178 | 399 | 77038751059 |
| 19 | 553 | 246 | 67 | 432 | 104505703249 |
| 20 | 661 | 147 | 127 | 503 | 192300806884 |
| 21 | 643 | 398 | 479 | 504 | 246215557249 |
| 22 | 721 | 423 | 317 | 587 | 366281426404 |
| 23 | 912 | 323 | 72 | 701 | 724451702259 |
| 24 | 931 | 679 | 529 | 813 | 1388951042164 |
|  | 1256 | 631 | 86 | 997 | 2964216377059 |
|  | 9631 | 1397 | 2429 | 7313 | 8615131736524564 |
|  | 16897 | 3348 | 16703 | 6064 | 81892097850563329 |
|  | 373283 | 141763 | 129541 | 286971 | 20627328664386373809604 |

**Table A3**. $a = 4$ (19 solutions, $n \leq 10^{11}$)

| $i$ | $A$ | $B$ | $C$ | $D$ | $n$ |
|---:|---:|---:|---:|---:|---:|
| 1 | 9 | 4 | 7 | 6 | 7585 |
| 2 | 47 | 3 | 33 | 31 | 4880005 |
| 3 | 61 | 32 | 19 | 46 | 18040145 |
| 4 | 107 | 73 | 101 | 77 | 244672565 |
| 5 | 137 | 14 | 103 | 88 | 352429025 |
| 6 | 173 | 1 | 163 | 83 | 895745045 |
| 7 | 111 | 161 | 63 | 163 | 2839400005 |
| 8 | 219 | 122 | 11 | 168 | 3186391345 |
| 9 | 237 | 58 | 61 | 168 | 3200222545 |
| 10 | 239 | 13 | 1 | 169 | 3262922885 |
| 11 | 303 | 59 | 241 | 189 | 8477361925 |
| 12 | 301 | 103 | 179 | 209 | 8658744725 |
| 13 | 339 | 14 | 259 | 216 | 13206989905 |
| 14 | 317 | 171 | 147 | 239 | 13518183445 |
| 15 | 243 | 254 | 163 | 264 | 20136041425 |
| 16 | 449 | 112 | 193 | 316 | 41272370945 |
| 17 | 469 | 127 | 331 | 311 | 49423420085 |
| 18 | 443 | 284 | 373 | 326 | 64535231345 |
| 19 | 555 | 94 | 165 | 392 | 95191700209 |
| 20 | 529 | 324 | 129 | 418 | 122390827585 |

**Table A4**. $a = 5$ (first 10 solutions, $n \leq 10^{12}$)

| $i$ | $A$ | $B$ | $C$ | $D$ | $n$ |
|---:|---:|---:|---:|---:|---:|
| 1 | 3 | 0 | 1 | 2 | 81 |
| 2 | 22 | 17 | 4 | 19 | 651861 |
| 3 | 83 | 52 | 19 | 64 | 84016401 |
| 4 | 147 | 63 | 71 | 101 | 545713686 |
| 5 | 190 | 11 | 40 | 127 | 1303283205 |
| 6 | 197 | 85 | 49 | 137 | 1767141606 |
| 7 | 417 | 117 | 19 | 281 | 31174327926 |
| 8 | 491 | 284 | 263 | 362 | 90647000241 |
| 9 | 563 | 114 | 409 | 348 | 101313827041 |
| 10 | 957 | 322 | 3 | 650 | 892531250081 |
| | | | | | |
| 11 | 1129 | 23 | 91 | 755 | 1624711078086 |
| 12 | 1396 | 233 | 1378 | 451 | 3812620279461 |
| 13 | 1613 | 338 | 423 | 1080 | 6834460387041 |
| | | | | | |
| | 4393 | 1091 | 3951 | 2283 | 379513978100406 |
| | 59201 | 10064 | 58879 | 15860 | 12334623789892762881 |
| | 64151 | 34620 | 51031 | 43152 | 24118655609148141601 |

**Table A5**. $a = 9$ (first 8 solutions, $n \le 10^{12}$)

| $i$ | $A$ | $B$ | $C$ | $D$ | $n$ |
|---|---|---|---|---|---|
| 1 | 29 | 9 | 11 | 17 | 766330 |
| 2 | 41 | 24 | 23 | 28 | 5811745 |
| 3 | 101 | 41 | 47 | 61 | 129492250 |
| 4 | 613 | 3 | 577 | 241 | 141202342090 |
| 5 | 625 | 77 | 85 | 361 | 152904267994 |
| 6 | 749 | 9 | 353 | 427 | 314722181050 |
| 7 | 830 | 329 | 250 | 503 | 580028236729 |
| 8 | 937 | 266 | 791 | 466 | 815887268785 |
| | | | | | |
| 9 | 896 | 489 | 148 | 599 | 1159123203625 |
| 10 | 589 | 738 | 499 | 742 | 2790087960865 |

**References**.

Hayashi, T. (1911), "On the Diophantine Equation $X^4 + Y^4 = Z^4 + T^4$," *Tohoku Mathematics Journal*, 1, 143-145

Izadi, F. and Baghalaghdam, M. (2017), "Is the Quartic Diophantine Equation $A^4 + hB^4 = C^4 + hD^4$ Solvable for any Integer $h$?", *arXiv:1701.02602v3*

Lander, L. J. (1968), "Geometric Aspects of Diophantine Equations Involving Equal Sums of Like Powers," *American Mathematics Monthly*, 75, 1061-1073

Brudno, S. (1969), "Some New Results on Equal Sums of Like Powers," *Mathematics of Computation*, 23, 877-880

Choudhry, A. (1998), "The Diophantine Equation $A^4 + 4B^4 = C^4 + 4D^4$ ", *Indian Journal of Pure and Applied Mathematics*, 29 (11), 1127-1128 (https://insa.nic.in/writereaddata/UpLoadedFiles/IJPAM/20005a05_1127.pdf)

Choudhry, A. (2016), "A Note on the Diophantine Equation $A^4 + hB^4 = C^4 + hD^4$ , *arXiv:1608.06220v1*

Piezas, T. On the Quartic Diophantine Equation $a^4 + nb^4 = c^4 + nd^4$ , https://sites.google.com/site/tpiezas/0021e

Roediger, P. A. (1972), "Notes on the Diophantine Equation $A^4 + aB^4 = C^4 + aD^4$ ," *(Unpublished, presented at the Eighteenth Conference of Army Mathematicians)*

Wroblewski, J. (Date Unknown), Magic Properties of Solutions of the Diophantine Equation $A^4 - B^4 = C^4 - D^4 = E^4 - F^4$ , http://www.math.uni.wroc.pl/~jwr/422/magic.pdf

Wroblewski, J. (Date Unknown), http://www.math.uni.wroc.pl/~jwr/eslp/tables.htm (422-10m.txt)

Zajta, A. J. (1983), "Solutions of the Diophantine Equation $A^4 + B^4 = C^4 + D^4$ ," *Mathematics of Computation*, v. 41, 635-659